\documentclass[opre,nonblindrev]{informs3} % current default for manuscript submission

\OneAndAHalfSpacedXI
\usepackage{amsmath, amsfonts,amssymb, latexsym,epsfig,color,rotating,paralist,
times,float,algorithm,verbatim,enumerate,algorithm,algorithmic, caption}
\usepackage{natbib}
\usepackage{endnotes}

 \def\newblock{\ }%
\usepackage[scriptsize]{subfigure}
\newcommand{\state}{x}
\renewcommand{\time}{k}
\newcommand{\statespace}{\mathbb{X}}
\newcommand{\statedim}{\textit{X}}
\newcommand{\actionspace}{\mathcal{A}}
\newcommand{\actiondim}{A}
\newcommand{\action}{a}
\newcommand{\tp}{P}
\newcommand{\obspace}{\mathbb{Y}}
\newcommand{\obsdim}{Y}
\newcommand{\oprob}{B}
\newcommand{\poly} {S}

\newcommand{\policy}{\mu}

\newcommand{\p}{\prime}
\newcommand{\E}                 {\mathbb{E}}
\newcommand{\prob}                 {\mathbb{P}}
\newcommand{\discount}{\rho}
\newcommand{\one}{\mathbf{1}}
\newcommand{\eye}{\textit{I}}
\newcommand{\belief}{\pi}
\newcommand{\Belief}{\Pi}
\newcommand{\tpe}{p}
\newcommand{\oprobe}{b}
\newcommand{\noise}{n}

\newcommand{\normal}{\mathcal{N}}

\newcommand{\copomat}{\Gamma}
\newcommand{\Real}{\mathbb{R}}
\newcommand{\filter}{\textit{T}}
\newcommand{\obs}{\textit{y}}
\newcommand{\uvalue}{\underline{V}}
\newcommand{\bvalue}{\overline{V}}
\newcommand{\valueaction}{Q}

\newcommand{\bvalueaction}{\overline{Q}}
\newcommand{\optvalue}{V}
\newcommand{\trid}{\Upsilon}
\newcommand{\filternorm}{\sigma}
\newcommand{\cost}{\textit{c}}

\newcommand{\gr}{\ge_{\text{r}}}
\newcommand{\gs}{\ge_{\text{s}}}
\newcommand{\gtp}{\ge_{\text{tp}}}
\newcommand{\gc} {\succeq}
\newcommand{\lc} {\preceq}

\newcommand{\lcost}{\underline{\textit{C}}}
\newcommand{\ucost}{\overline{\textit{C}}}
\newcommand {\uf} {\textit{f}}%{\underline{\textit{f}}\xspace}
\newcommand {\of} {\textit{g}}%{\overline{\textit{f}}\xspace}

\newcommand {\setoverlap} {{\Belief_{O}}}
\newcommand {\optpolicy} {\mu^*}
\newcommand {\policyu} {\overline{\mu}}
\newcommand {\policyl} {\underline{\mu}}
\newcommand {\diag} {\text{diag}}

\newcommand{\basisvec} {\textit{e}}

\newcommand {\percentloss} {\epsilon}
\newcommand {\approxpolicy} {\tilde{\mu}}
\usepackage{epsfig}
\usepackage{xspace}
\usepackage[margin=1in]{geometry}
\usepackage{amsmath}
\usepackage{amsfonts}
\newtheorem{theorem}{Theorem}
\newtheorem{lemma}[theorem]{Lemma}

\begin{document}

\RUNTITLE{Myopic Bounds for POMDPs}
\RUNAUTHOR{Krishnamurthy}
%\title{Myopic Upper and Lower Bounds for the Optimal Policy of Partially Observed Markov Decision Process}
\TITLE{Myopic Bounds for  Optimal Policy of POMDPs: An extension of Lovejoy's structural results}

\ARTICLEAUTHORS{%
\AUTHOR{Vikram Krishnamurthy and Udit Pareek}
\AFF{Department of Electrical and Computer Engineering, University of British Columbia, Vancouver, V6T 1Z4, Canada. vikramk@ece.ubc.ca}}
\MANUSCRIPTNO{published in Vol.63, Jan. 2015}

\ABSTRACT{
This paper provides a relaxation of the sufficient conditions, and also an  extension of the structural results for  Partially Observed Markov Decision Processes (POMDPs) given in \cite{Lov87}. Sufficient conditions are provided so that  the optimal policy  can be upper and lower bounded  by judiciously chosen myopic policies. These myopic policy bounds are constructed to maximize the volume of belief states where they coincide with the optimal policy. 
Numerical examples illustrate these myopic bounds for both continuous and discrete observation sets.}


\maketitle

\section{Introduction} \label{sec:intro}
POMDPs have received much attention due to their applications in diverse areas such as scheduling in sensor networks and wireless communications (see \cite{Kri11}, \cite{KD07} and references therein) and  artificial intelligence \citep{KLC98}.  Even though, for  finite observation alphabet sets, and finite horizon, the optimal policy of a POMDP can be computed via stochastic dynamic programming, such problems are P-SPACE hard \citep{PT87}.

The seminal papers \cite{Lov87,Rie91,RZ94} give sufficient conditions such that the optimal policy of a POMDP can be lower bounded by a myopic policy. Unfortunately, despite the  enormous usefulness of such
a result, the sufficient conditions given in \cite{Lov87} and \cite{Rie91} for this  result to hold are not useful - it is impossible to generate non-trivial examples that satisfy the conditions (c), (e), (f) of \cite[Proposition 2]{Lov87} and condition (i) of \cite[Theorem 5.6]{Rie91}. In this paper, we provide a fix to these sufficient conditions so that  the results of \cite{Lov87,Rie91} hold for constructing a myopic policy that lower bounds the optimal policy. Then, for infinite horizon discounted cost POMDPs, we show how this idea of constructing a  lower bound myopic policy can be extended to constructing an {\em upper} bound myopic policy. More specifically, for belief state  $\belief$, we present sufficient conditions under which the optimal policy, denoted by $\optpolicy(\belief)$, of a given POMDP can be upper and lower bounded by myopic policies denoted by $\policyu(\belief)$ and $\policyl(\belief)$, respectively, i.e. $\policyl(\belief) \le \optpolicy(\belief) \le \policyu(\belief)$ for all $ \belief \in \Belief$. Here $\Belief$ denotes the set of belief states of a POMDP. Interestingly, these judiciously chosen myopic policies are independent of the actual values of the observation probabilities (providing they satisfy a sufficient condition) which makes the structural results  applicable to both discrete and continuous observations.
Finally, we construct the myopic policies, $\policyu(\belief)$ and $\policyl(\belief)$, to maximize the volume of the belief space where they coincide with the optimal policy $\mu^*(\pi)$.

Numerical examples are presented to illustrate the performance of  these myopic policies.
To quantify how well the myopic policies perform we use two parameters:  the volume of the belief space where the myopic policies coincide with the optimal policy, and
an upper bound to the average percentage loss in optimality due to following this optimized myopic policy.

\section {The Partially Observed Markov Decision Process} \label{sec:pomdp}
Consider  a discrete time, infinite horizon discounted cost POMDP. A   discrete time Markov chain  evolves on the  state space $\statespace = \{1,2,\ldots, \statedim\}$. Denote the
action space  as $\actionspace = \{1,2,\ldots,\actiondim\}$ and observation space as $\obspace$. For discrete-valued observations $\obspace = \{1,2,\ldots,\obsdim\}$ and for continuous observations $\obspace \subset \Real$.

Let
$\Belief = \left\{\belief: \belief(i) \in [0,1], \sum_{i=1}^\statedim \belief(i) = 1 \right\}$ denote the belief space of $\statedim$-dimensional probability vectors.  For stationary policy  $\policy: \Belief \rightarrow \actionspace$,
 initial belief  $\belief_0\in \Belief$,  discount factor $\discount \in [0,1)$, define the  discounted cost:
\begin{align}\label{eq:discountedcost}
J_{\policy}(\belief_0) = \E\left\{\sum_{\time=1}^{\infty} \discount ^{\time-1} \cost_{\policy(\belief_\time)}^\p\belief_\time\right\}.
\end{align}
%that minimizes (\ref{eq:discountedcost}).
Here $\cost_\action = [\cost(1,\action),\ldots,\cost(\statedim,\action)]^\p$, $a\in \actionspace$ is the cost vector for each action, and the belief state evolves as
$\belief_{k} = \filter(\belief_{k-1},\obs_k,\action_k)$ where
\begin{align}  \filter\left(\belief,\obs,\action\right) = \cfrac{\oprob^{\action}_{\obs} \tp'_{\action}\belief}{\filternorm\left(\belief,\obs,\action\right)} , \quad
\filternorm\left(\belief,\obs,\action\right) = \one_{\statedim}'\oprob^{\action}_{\obs} \tp'_{\action}\belief, \quad
\oprob^{\action}_{\obs} = \diag\{\oprobe^{\action}_{1,\obs},\cdots,\oprobe^{\action}_{\statedim,\obs}\}. \label{eq:information_state}
\end{align}
Here  $\one_{\statedim}$ represents a $\statedim$-dimensional vector of ones,
$ \tp_\action = \left[\tpe^\action_{ij}\right]_{\statedim\times\statedim}$
$ \tpe_{ij}^{\action} = \prob(\state_{\time+1} = j | \state_\time = i, \action_\time =a )$ denote the transition probabilities,
 $\oprobe_{\state\obs}^\action = \prob(\obs_{\time+1} = \obs| \state_{\time+1} = \state, \action_{\time} = \action)$ when $\obspace$ is finite,
 or  $\oprobe^\action_{\state\obs}$ is the conditional probability density function when $\obspace \subset \Real$.

The aim is to compute the optimal  stationary policy $\optpolicy:\Belief \rightarrow \actionspace$ such that
$J_{\optpolicy}(\belief_0) \leq J_{\policy}(\belief_0)$ for all $\belief_0 \in \Belief$.
Obtaining the optimal policy  $\optpolicy$ is equivalent to solving
 Bellman's  dynamic programming equation:
$ \optpolicy(\belief) =  \underset{\action \in \actionspace}\argmin ~\valueaction(\belief,\action)$, $J_{\optpolicy}(\belief_0) = \optvalue(\belief_0)$, where
\begin{equation}
\optvalue(\belief)  = \underset{\action \in \actionspace}\min ~\valueaction(\belief,\action), \quad
  \valueaction(\belief,\action) =  ~\cost_\action^\prime\belief + \discount\sum_{\obs \in \obsdim} \optvalue\left(\filter\left(\belief,\obs,\action\right)\right)\filternorm \left(\belief,\obs,\action\right). \label{eq:bellman}
\end{equation}

Since  $\Belief$ is continuum, Bellman's equation \eqref{eq:bellman} does not translate into practical solution methodologies as $\optvalue(\belief)$ needs to be evaluated at each $\belief \in \Belief$. This motivates the construction of  judicious upper and lower bounds, denoted by ${\policyu(\belief)}$ and ${\policyl(\belief)}$ respectively, to the optimal policy $\optpolicy(\belief)$. For belief states $\belief$ where ${\policyu(\belief)}={\policyl(\belief)}$, the optimal policy $\optpolicy(\belief)$ is completely determined.

\section {Myopic Bounds to the Optimal Policy} \label{sec:exist_bounds}
\subsection{Assumptions}\label{subsec:assumptions}

\begin{enumerate}[\bf{(A}1)]
\item\label{it:increasing_cost} There exists a $\of \in \Real^{\statedim}$ such that $\ucost_\action \equiv \cost_{\action} + \left(\eye - \discount\tp_\action\right)\of$ is strictly increasing in $\state \in \statespace,\forall \action \in \actionspace$.
\item\label{it:decreasing_cost} There exists a $\uf \in \Real^{\statedim}$ such that $\lcost_\action \equiv \cost_{\action} + \left(\eye - \discount\tp_\action\right)\uf$ is strictly decreasing in $\state \in \statespace,\forall \action \in \actionspace$.
\item\label{it:TP2_gen} $\tp_\action$ and $\oprob_\action$, $\action \in \actionspace$ are totally positive of order 2 (TP2), that is, all second-order minors are nonnegative.
\item\label{it:Copo_gen}
$ \gamma^{j,\action,\obs}_{mn} + \gamma^{j,\action,\obs}_{nm} \ge 0~\forall m,n,j,\action,\obs ~\text{where}~ \gamma^{j,\action,\obs}_{mn} =  \oprobe^\action_{j,\obs}\oprobe^{\action+1}_{j+1,\obs}\tpe^{\action}_{m,j}\tpe^{\action+1}_{n,j+1} - \oprobe^\action_{j+1,\obs}\oprobe^{\action+1}_{j,\obs}\tpe^{\action}_{m,j+1}\tpe^{\action+1}_{n,j}$.
\item \label{it:norm_ac_inc}$\sum_{\obs \le \bar{\obs}}\sum_{j \in \statespace}\left[\tpe^{\action}_{i,j}\oprobe^{\action}_{j,\obs} - \tpe^{\action+1}_{i,j}\oprobe^{\action+1}_{j,\obs}\right] \le 0 \quad \forall i \in \statespace, \forall \bar{\obs} \in \obspace$
\end{enumerate}

\subsubsection*{Discussion}
If  the elements of $\cost_\action$ are strictly increasing then \textbf{(A\ref{it:increasing_cost})} holds trivially. Similarly, if the elements of $\cost_\action$ are decreasing then \textbf{(A\ref{it:decreasing_cost})} holds and coincides with Assumption (b) in \cite[Proposition 2]{Lov87}.

 \textbf{(A\ref{it:increasing_cost})} and \textbf{(A\ref{it:decreasing_cost})} are easily verified by checking the feasibility of the following linear programs:
    \begin{align}\label{eq:lp_a1_a2}
    LP1:&  \underset{\of \in \poly_{\of}}{\min}~ \one_{\statedim}^\prime\of,~LP2: \underset{\uf \in \poly_{\uf}}{\min} ~\one_{\statedim}^\prime\uf. \\
\label{eq:polytopes_u}
    \poly_{\of} &= \left\{ \of:  \ucost_\action^\prime\basisvec_i \le \ucost_\action^\prime\basisvec_{i+1} ~~\forall \action \in \actionspace, i \in \statespace \right \}  \\
  \label{eq:polytopes_l}
\poly_{\uf} &= \Big\{ \uf:  \lcost_\action^\prime\basisvec_i \ge \lcost_\action^\prime\basisvec_{i+1} ~~\forall \action \in \actionspace, i \in \statespace \Big \}
\end{align}
where $\basisvec_i$ is the unit  $\statedim$-dimensional vector with 1 at the $i$th position.

 \textbf{(A\ref{it:TP2_gen})}  is equivalent to saying that the rows of $\tp_\action$ and $\oprob_\action$ are monotone likelihood ratio (MLR) increasing.
  (MLR dominance is defined in the appendix).
 Numerous examples of  TP2 matrices satisfying \textbf{(A\ref{it:TP2_gen})} can be found in \cite{KR80}. Examples of TP2 observation kernels include Gaussian, Exponential, Binomial and  Poisson distributions. Examples of discrete observation distributions include binary erasure channels and binary symmetric channels with error probability less than 0.5.

\textbf{(A\ref{it:Copo_gen})}  implies that   the belief due to action $\action+1$  MLR dominates the  belief due to action $\action$, i.e., in the terminology of
  \cite{Mil81}, $\action+1$ yields a more ''favorable outcome"  than  $\action$. 
For POMDPs with $\oprob_\action = \oprob~\forall \action \in \actionspace$, \textbf{(A\ref{it:Copo_gen})} trivially holds for TP2 transition matrices $\tp_\action$ and $\tp_{\action'}$, $\action >\action'$  if all rows of $\tp_\action$ MLR dominate the last row of $\tp_{\action'}$.

  \textbf{(A\ref{it:Copo_gen})} and \textbf{(A\ref{it:norm_ac_inc})} are a relaxed version of Assumptions (c), (e), (f) of \cite[Proposition 2]{Lov87} and Assumption (i) of \cite[Theorem 5.6]{Rie91}.
  In particular, the assumptions (c), (e), (f) of \cite{Lov87} require that $\tp_{\action+1}$ $ \gtp \tp_{\action}$ and $\oprob_{\action+1} \gtp \oprob_{\action}$, where $\gtp$ (TP2 stochastic ordering) is defined in \cite{MS02}, which is impossible for stochastic matrices, unless $\tp_\action =  \tp_{\action+1}$, $\oprob_\action =\oprob_{\action+1}$ or the matrices $\tp_\action, \oprob_\action$ are rank 1 for all $\action$ meaning that the observations are non-informative.

  Assumptions (c) and (f) of \cite[Proposition 2]{Lov87} are required to ensure that the posterior $\filter(\belief,\obs,\action)$ \eqref{eq:information_state} is MLR increasing in $\action$. A necessary and sufficient condition to ensure the monotonicity of $\filter(\belief,\obs,\action)$ is that the matrices $\copomat^{j,\action,\obs}$, defined below, are copositive \cite{BD09} on $\Belief$. That is,
\begin{align}\label{eq:copositive}
\begin{aligned}
&\belief^\prime\copomat^{j,\action,\obs}\belief \ge 0, \forall \belief \in \Belief, \forall j, \action, \obs~\text{where}\\
&\copomat^{j,\action,\obs} = \cfrac{1}{2}\left[\gamma^{j,a,y}_{mn} + \gamma^{j,a,y}_{nm}\right]_{\statedim\times\statedim},\gamma^{j,a,y}_{mn} = \oprobe^\action_{j,\obs}\oprobe^{\action+1}_{j+1,\obs}\tpe^{\action}_{m,j}\tpe^{\action+1}_{n,j+1} - \oprobe^\action_{j+1,\obs}\oprobe^{\action+1}_{j,\obs}\tpe^{\action}_{m,j+1}\tpe^{\action+1}_{n,j}.
\end{aligned}
\end{align}
In general, the problem of verifying the copositivity of a matrix is NP-complete. Assumption \textbf{(A\ref{it:Copo_gen})} is a simpler but more restrictive sufficient condition to ensure that $\copomat^{j,\action,\obs}$ \eqref{eq:copositive} is copositive.

\subsection{Construction of Myopic Upper and Lower Bounds}\label{subsec:subsec_asmp_result}
We are interested in   myopic policies of the form $\underset{\action \in \actionspace}{\argmin}~ C_\action^\prime\pi$ where cost vectors $C_\action$
are constructed so that when applied to Bellman's equation \eqref{eq:bellman}, they leave the optimal policy $\mu^*(\pi)$ unchanged. This is for several reasons: First, similar to \cite{Lov87}, \cite{Rie91} it allows us to construct useful myopic policies that provide provable upper and lower bounds to the optimal policy.  Second, these myopic policies can be straightforwardly extended to 2-stage or multi-stage myopic costs. Third, such a choice  precludes choosing useless myopic bounds such
as $\policyu(\belief) = \actiondim$ for all $ \belief \in \Belief$.

Accordingly, for any two vectors $\of$ and $\uf \in \Real^\statedim$, define the myopic policies associated with the transformed costs $\ucost_{\action}$ and $\lcost_{\action}$ as follows:
\begin{align}\label{eq:new_cost_u}
\policyu(\belief) & \equiv \underset{\action \in \actionspace}\argmin ~\ucost_{\action}'\belief,\quad \text{ where } ~\ucost_{\action} = \cost_\action + \left( \eye - \discount \tp_\action \right)\of \\
\label{eq:new_cost_l}
\policyl(\belief)& \equiv \underset{\action \in \actionspace}\argmin ~\lcost_{\action}'\belief,\quad \text{ where } ~\lcost_{\action} = \cost_\action + \left( \eye - \discount \tp_\action \right)\uf.
\end{align}
It is easily seen that Bellman's equation \eqref{eq:bellman} applied to optimize the objective \eqref{eq:discountedcost} with transformed costs $\ucost_\action$ and $\lcost_\action$ yields the same optimal strategy $\optpolicy(\belief)$ as the Bellman's equation with original costs $\cost_\action$. The corresponding value functions are $\bvalue(\belief) \equiv \optvalue(\belief) + \of^\prime\belief$ and $\uvalue(\belief) \equiv \optvalue(\belief) + \uf^\prime\belief$. The following  main result is proved  in the Appendix.
\begin{theorem}
\label{th:theorem1}
Consider a POMDP $\left(\statespace, \actionspace, \obspace, \tp_\action, \oprob_\action, \cost, \discount\right)$ and assume \textbf{(A\ref{it:increasing_cost})}-\textbf{(A\ref{it:norm_ac_inc})} holds. Then the myopic policies, $\policyu(\belief)$ and $\policyl(\belief)$, defined in \eqref{eq:new_cost_u}, \eqref{eq:new_cost_l} satisfy:
$\policyl(\belief) \le \optpolicy(\belief) \le \policyu(\belief)$ for all  $\belief \in \Belief$.
\end{theorem}

\section {Optimizing the Myopic Policy Bounds to Match the Optimal Policy} \label{sec:optimize_bounds}
The aim of this section is to determine vectors $\of$ and $\uf$, in \eqref{eq:polytopes_u} and \eqref{eq:polytopes_l}, that maximize the volume of the simplex where the myopic upper and lower policy bounds, specified by \eqref{eq:new_cost_u} and \eqref{eq:new_cost_l}, coincide with the optimal policy.
That is, we wish to maximize the volume of the `overlapping region'
\begin{align}\label{eq:def_overlap}
\begin{aligned}
\text{vol}\left(\setoverlap\right),~\text{where}~ \setoverlap \equiv \{\belief: \policyu(\belief) = \policyl(\belief) = \optpolicy(\belief)\}.
\end{aligned}
\end{align}
Notice that the myopic policies $\policyu$ and $ \policyl$ defined in (\ref{eq:new_cost_u}), (\ref{eq:new_cost_l})  do not depend on the observation probabilities $\oprob_\action$ and so
neither does $\text{vol}\left(\setoverlap\right)$.  So  $\policyu$ and  $ \policyl$ can be chosen to maximize $\text{vol}\left(\setoverlap\right)$ independent
of $\oprob_\action$ and therefore work for  discrete and continuous observation spaces.  Of course, the proof of Theorem  \ref{th:theorem1} requires conditions on $\oprob_\action$.

\subsection {Optimized Myopic Policy for Two Actions} \label{subsec:optimize_two_actions}
For a two action POMDP, obviously for a belief $\belief$, if $\policyu(\belief)=1$ then $\optpolicy(\belief) = 1$. Similarly, if  $\policyl(\belief)=2$, then $ \optpolicy(\belief) = 2$.
 Denote the set of beliefs (convex polytopes) where $\policyu(\belief) = \optpolicy(\belief) = 1$ and $\policyl(\belief)= \optpolicy(\belief) = 2$ as
\begin{align}\label{eq:overlapping_regions}
\begin{aligned}
\Belief^{\of}_{1} &= \left\{ \belief:  \ucost_1^\prime\belief \le \ucost_2'\belief \right\}= \left\{ \belief:  (\cost_1-\cost_2 - \discount(\tp_1 - \tp_2)\of)'\belief  \le 0\right\},\\
\Belief^{\uf}_{2} &= \left\{ \belief:  \lcost_2'\belief \le \lcost_1'\belief \right\}= \left\{ \belief:  (\cost_1-\cost_2 - \discount(\tp_1 - \tp_2)\uf)'\belief  \ge 0\right\}.
\end{aligned}
\end{align}
Clearly  $\setoverlap = \Belief^{\of}_{1} \cup \Belief^{\uf}_{2}$. Our goal is to find $\of^* \in \poly_{\of}$ and $\uf^* \in \poly_{\uf}$ such that $\text{vol}\left(\setoverlap\right)$ is maximized.
\begin{theorem}
\label{th:theorem_opt_bounds}
Assume that there exists two fixed $\statedim$-dimensional vectors $\of^*~\text{and}~\uf^*$ such that
\begin{align}\label{eq:optimal_overlap}
\begin{aligned}
&(\tp_2 - \tp_1)\of^* \lc (\tp_2 - \tp_1)\of , \; \forall \of \in \poly_{\of}\\
&(\tp_1 - \tp_2)\uf^* \lc (\tp_1 - \tp_2)\uf, \; \;\forall \uf \in \poly_{\uf}
\end{aligned}
\end{align}
where for $\statedim$-dimensional vectors $a$ and $b$, $a \lc b \Rightarrow \left[a_1 \le b_1,\cdots, a_{\statedim} \le b_{\statedim}\right]$,
and $ \poly_{\of},  \poly_{\uf}$ are defined in (\ref{eq:polytopes_u}), (\ref{eq:polytopes_l}), respectively.
If  the myopic policies $\policyu$ and $\policyl$ are constructed using $\of^*$ and $\uf^*$, then $\text{vol}(\setoverlap)$ is maximized.
\end{theorem}

Theorem  \ref{th:theorem_opt_bounds} asserts that  myopic policies $\policyu$ and $\policyl$ characterized by two
fixed vectors $\of^*~\text{and}~\uf^*$
maximize  $\text{vol}(\setoverlap)$ over the entire belief space $\Belief$.
The existence and computation of these policies characterized by
 $\of^* \in \poly_{\of}$ and $\uf^* \in \poly_\uf$ are determined by  Algorithm~\ref{alg1}.
Algorithm \ref{alg1} solves $\statedim$ number of LPs to obtain $\of^*$. If no $\of^* \in \poly_{\of}$ satisfying \eqref{eq:optimal_overlap} exists, then Algorithm \ref{alg1} will terminate with no solution. The procedure for computing $\uf^*$ is similar. \begin{algorithm}
\caption{Compute $\of^*$}
\label{alg1}
\begin{algorithmic}[1]
%\STATE $\of^*\gets$ ComputeMinimal($\tp^1,\tp^2,\discount,\cost$)
 \FORALL{$i \in \statedim$}
  \STATE \label{op0}$\alpha_i \gets \underset{\of \in \poly_{\of}}\min ~\basisvec_i'(\tp_2 - \tp_1)\of$
 \ENDFOR
 \STATE $\of^* \in \poly_{\of^*}, \poly_{\of^*} \equiv \left\{\of^*:\of^*\in \poly_{\of}, \basisvec_i'(\tp_2 - \tp^1)\of^* = \alpha_i, i = 1,\cdots,\statedim\right\}  $\label{op1}
 \STATE $\policyu(\belief) = \underset{\action \in \{1,2\}}{\argmin}~\belief'\ucost^{*}_{\action} ~\forall \belief \in \Belief$, where $\ucost^{*}_{\action} = \cost_\action + \left( \eye - \discount \tp_\action \right)\of^*$
 \STATE $\policyu(\belief) = \optpolicy(\belief) = 1, \forall \belief \in \Belief^{\of^*}_{1}$.
\end{algorithmic}
\end{algorithm}

\subsection {Optimizing Myopic Policies for more than 2 actions} \label{subsec:optimize_gr_2}Unlike Theorem \ref{th:theorem_opt_bounds},
for the case  $\actiondim > 2$, we are unable to show that a single fixed choice of $\policyu$ and $\policyl$  maximizes $\text{vol}(\setoverlap)$. Instead at each time $\time$,  $\policyu$ and $\policyl$ are optimized
depending on the belief state $\belief_\time$.
Suppose at time $\time$, given observation $\obs_\time$, the belief state, $\belief_\time$, is computed by using  \eqref{eq:information_state}. For this belief state $\belief_\time$, the aim is to compute $\of^* \in \poly_{\of}$ \eqref{eq:polytopes_u} and $\uf^* \in \poly_{\uf}$ \eqref{eq:polytopes_l} such that the difference between myopic policy bounds, $\policyu(\belief_\time) - \policyl(\belief_\time)$, is minimized. That is,
\begin{align}\label{eq:LP_min_dist}
\begin{aligned}
    \left(\of^*,\uf^*\right) = \underset{\of \in \poly_{\of},~\uf \in \poly_\uf}{\argmin} ~\policyu(\belief_\time) - \policyl(\belief_\time).
\end{aligned}
\end{align}
\eqref{eq:LP_min_dist} can be decomposed into following two optimization problems,
\begin{align}\label{eq:LP_min_ub}
\begin{aligned}
    \of^* = \underset{\of \in \poly_{\of}}{\argmin}~ \policyu(\belief_\time), ~\uf^* = \underset{\uf \in \poly_{\uf}}{\argmax}~ \policyl(\belief_\time).
\end{aligned}
\end{align}
If assumptions \textbf{(A\ref{it:increasing_cost})} and \textbf{(A\ref{it:decreasing_cost})} hold, then the optimizations in \eqref{eq:LP_min_ub} are feasible.
Then $\policyu(\belief_\time)$ in \eqref{eq:new_cost_u} and $\of^*$, in \eqref{eq:LP_min_ub} can be computed as follows:  Starting with $\policyu(\belief_\time) = 1$, successively solve a maximum of $\actiondim$ feasibility LPs, where the $i$th LP searches for a feasible $\of \in \poly_\of$ in \eqref{eq:polytopes_u} so that the myopic upper bound yields action $i$, i.e. $\policyu(\belief_\time) = i$. The $i$th feasibility LP can be written as
 \begin{align}\label{eq:LP_UB_iterative}
\begin{aligned}
    &~\underset{\of \in \poly_\of}{\min}~\one_\statedim^\prime\of\\
    &s.t.,~~~\ucost_i^\prime\belief_\time \le \ucost_\action^\prime\belief_\time~\forall \action \in \actionspace, \action \ne i
\end{aligned}
\end{align}
The smallest $i$, for which \eqref{eq:LP_UB_iterative} is feasible, yields the solution $\left(\of^*, \policyu(\belief_\time)=i\right)$ of the optimization in \eqref{eq:LP_min_ub}. The above procedure is straightforwardly modified to obtain $\uf^*$ and the lower bound $\policyl(\belief_\time)$ \eqref{eq:new_cost_l}.

\section {Numerical Examples} \label{sec:numerical_examples}
Recall that on the set $\setoverlap$ \eqref{eq:def_overlap}, the upper and lower myopic bounds coincide with the optimal policy $\optpolicy(\belief)$. What is the performance loss outside the set $\setoverlap$? To quantify this, define the  policy
\begin{align}\nonumber
\begin{aligned}
\approxpolicy(\belief) = \begin{cases} \optpolicy(\belief) &\forall \belief \in \setoverlap\\
\text{arbitrary action (e.g. 1)} &\forall \belief \not\in \setoverlap\end{cases}
\end{aligned}
\end{align}
Let  $J_{\approxpolicy}(\belief_0)$ denote the discounted cost associated with $\approxpolicy(\belief_0)$. Also denote
\begin{align}\nonumber
\begin{aligned}
\tilde{J}_{\optpolicy}(\belief_0) = \E\left\{\sum_{k=1}^{\infty} \discount^{k-1} \tilde{\cost}_{\optpolicy(\belief_\time)}^\p \belief_k\right\},~
\text{where}, \tilde{\cost}_{\optpolicy(\belief)} = \begin{cases} {\cost}_{\optpolicy(\belief)}& \belief \in \setoverlap\\
\left[\underset{\action \in \actionspace}{\min}~\cost(1,\action),\cdots,\underset{\action \in \actionspace}{\min}~\cost(\statedim,\action)\right]^\p& \belief \not\in \setoverlap\end{cases}
\end{aligned}
\end{align}
Clearly an upper bound for the  percentage loss in optimality due to using  policy $\approxpolicy$ instead of  optimal policy $\optpolicy$ is
\begin{align}\label{eq:percentloss}
\begin{aligned}
\percentloss = \cfrac{J_{\approxpolicy}(\belief_0) - \tilde{J}_{\optpolicy}(\belief_0)}{\tilde{J}_{\optpolicy}(\belief_0)}.
\end{aligned}
\end{align}
In the numerical examples below,  to evaluate $\percentloss$, 1000 Monte-Carlo simulations were run to estimate the discounted costs $J_{\approxpolicy}(\belief_0)$ and $\tilde{J}_{\optpolicy}(\belief_0)$ over a horizon of 100 time units.
%The expectation in \eqref{eq:percentloss} is evaluated by Monte-Carlo simulations over a set of 100 uniformly sampled priors $\belief_0$ from the simplex $\Belief$.
The parameters $\percentloss$ and $\text{vol}\left(\setoverlap\right)$ are used to evaluate the performance of the optimized myopic policy bounds constructed according to Sec. \ref{sec:optimize_bounds}. Note that  $\percentloss$ depends on the choice of observation distribution $\oprob$, unlike $\text{vol}\left(\setoverlap\right)$, see discussion below (\ref{eq:def_overlap}) and also Example 2 below.

\textit{Example 1. Sampling and Measurement Control with Two Actions}: In this problem \citep{Kri13}, at every decision epoch, the decision maker has the option of either recording a noisy observation (of a Markov chain) instantly (action $\action = 2$) or waiting for one time unit and then recording an  observation using a better sensor (action $\action = 1$). Should one record observations more frequently and less accurately or more accurately but less frequently?

We chose
 $\statedim = 3$, $\actiondim = 2$ and $\obsdim = 3$. Both  transition and observation probabilities are action dependent (parameters specified in the Appendix). The percentage loss in optimality is evaluated by simulation for different values of the discount factor $\discount$. Table \ref{tb:table1a} displays $\text{vol}\left(\setoverlap\right)$, $\percentloss_1$ and $\percentloss_2$. For each $\discount$, $\percentloss_1$ is obtained by assuming $\belief_0 = \basisvec_3$ (myopic bounds overlap at $\basisvec_3$) and $\percentloss_2$ is obtained by uniformly sampling $\belief_0 \notin \setoverlap$. Observe that $\text{vol}\left(\setoverlap\right)$ is large and $\percentloss_1$, $\percentloss_2$ are  small,  which indicates the usefulness of the proposed myopic policies.

\begin{table}[h!]
\caption{Performance of optimized myopic policies versus discount factor $\discount$ for four numerical examples.  The performance metrics $\text{vol}\left(\setoverlap\right)$ 
and $\percentloss$ are 
defined in (\ref{eq:def_overlap}) and (\ref{eq:percentloss}).}
\centering
%\subtable{
%\centering
%\begin{tabular}{c}\hline
%$\discount$ \\ \hline
% 0.5\\
% 0.6\\
% 0.7\\
% 0.8\\
% 0.9\\\hline
%\end{tabular}
%%\label{tb:table1_discount}
%}
\subtable[Example 1]{
\centering
\begin{tabular}{c|ccc}\hline
$\discount$ & \text{vol}$\left(\setoverlap\right)$ & $\percentloss_1$ & $\percentloss_2$\\ \hline
0.4& $95.3\%$  &  $0.30\% $ & $16.6\%$  \\
0.5& $94.2\%$  &  $0.61\%$ & $13.9\%$ \\
0.6& $92.4\%$  &  $1.56\%$ & $11.8\%$ \\
0.7& $90.2\%$  &  $1.63\%$ & $9.1\%$ \\
0.8& $87.4\%$  &  $1.44\%$ & $6.3\%$ \\
0.9& $84.1\%$  &  $1.00\%   $ & $3.2\%$    \\ \hline
\end{tabular}
\label{tb:table1a}
}
\subtable[Example 2]{
\centering
\begin{tabular}{ccccc}\hline
\text{vol}$\left(\setoverlap\right)$ & $\percentloss^d_1$&$\percentloss^d_2$&$\percentloss^c_1$&$\percentloss^c_2$  \\ \hline
$64.27\%$ & 7.73\%   & 12.88\% & 6.92\%   & 454.31\%  \\
$55.27\%$ & 8.58\%   & 12.36\% & 8.99\%   & 298.51\%  \\
$46.97\%$ & 8.97\%   & 11.91\% & 12.4\%   & 205.50\%  \\
$39.87\%$ & 8.93\%   & 11.26\% & 14.4\%  & 136.31\%  \\
$34.51\%$ & 10.9\%  & 12.49\% & 17.7\%  & 88.19\%  \\
$29.62\%$ & 11.2\%  & 12.24\% & 20.5\%  & 52.16\%  \\ \hline
\end{tabular}
\label{tb:table1c}
}
\subtable[Example 3]{
\centering
\begin{tabular}{ccc}\hline
\text{vol}$\left(\setoverlap\right)$ & $\percentloss_1$ & $\percentloss_2$\\ \hline
 $61.4\%$ & 2.5\% & 10.1\%	\\
 $56.2\%$ & 2.3\% & 6.9\%	\\
 $47.8\%$ & 1.7\% & 4.9\%	\\
 $40.7\%$ & 1.4\% & 3.5\%	\\
 $34.7\%$ & 1.1\% & 2.3\%	\\
 $31.8\%$ & 0.7\% & 1.4\%	\\ \hline
\end{tabular}
\label{tb:table1d}
}
\subtable[Example 4]{
\centering
\begin{tabular}{cccccc}\hline
$\overline{\text{vol}}\left(\setoverlap\right)$ & $\underline{\text{vol}}\left(\setoverlap\right)$ & $\overline{\percentloss}_1$ & $\underline{\percentloss}_1$ & $\overline{\percentloss}_2$ & $\underline{\percentloss}_2$\\ \hline
$98.9\%$ & $84.5\%$ & 0.10\% & 6.17\%	& 1.45\% & 1.71\%\\
$98.6\%$ & $80.0\%$ & 0.18\% & 7.75\%	& 1.22\% & 1.50\%\\
$98.4\%$ & $75.0\%$ & 0.23\% & 11.62\%	& 1.00\% & 1.31\%\\
$98.1\%$ & $68.9\%$ & 0.26\% & 14.82\%	& 0.75\% & 1.10\%\\
$97.8\%$ & $61.5\%$ & 0.27\% & 19.74\%	& 0.51\% & 0.89\%\\
$97.6\%$ & $52.8\%$ & 0.25\% & 24.08\%	& 0.26\% & 0.61\%\\ \hline
\end{tabular}
\label{tb:table1f}
}
\label{tb:table1}
\end{table}

\textit{Example 2. 10-state POMDP}: Consider a POMDP with $\statedim=10$, $\actiondim = 2$. Consider two sub-examples: the first with discrete observations  $\obsdim = 10$ (parameters  in Appendix), the second with  continuous  observations obtained using the additive Gaussian noise model, i.e. $\obs_\time = \state_\time + \noise_\time$ where $\noise_\time \thicksim \normal(0,1)$. The  percentage loss in optimality  is evaluated by simulation for these two sub examples and denoted by  $\percentloss^d_1, \percentloss^d_2$ (discrete observations) and $\percentloss^c_1, \percentloss^c_2$ (Gaussian observations)
in Table \ref{tb:table1c}.

$\percentloss^d_1$ and $\percentloss^c_1$ are obtained by assuming $\belief_0 = \basisvec_5$ (myopic bounds overlap at $\basisvec_5$). $\percentloss^d_2$ and $\percentloss^c_2$ are obtained by sampling $\belief_0 \notin \setoverlap$. Observe  from Table \ref{tb:table1c} that $\text{vol}\left(\setoverlap\right)$  decreases with $\discount$. However, the values of $\percentloss^d_1$ and $\percentloss^c_1$ are  small for all values of $\discount$ indicating  the usefulness of the  myopic bounds when the prior $\belief_0 \in \setoverlap$.

\textit{Example 3. 8-state and 8-action POMDP}: Consider a  POMDP with $\statedim=8$, $\actiondim = 8$ and $\obsdim = 8$ (parameters in  Appendix). Table \ref{tb:table1d} displays $\text{vol}\left(\setoverlap\right)$, $\percentloss_1$ and $\percentloss_2$. For each $\discount$, $\percentloss_1$ is obtained by assuming $\belief_0 = \basisvec_1$ (myopic bounds overlap at $\basisvec_1$) and $\percentloss_2$ is obtained by uniformly sampling $\belief_0 \notin \setoverlap$. The results indicate that the myopic policy bounds are still useful for some values of $\discount$.

\textit{Example 4.  Myopic Bounds versus Transition Matrix}:  The aim here is to illustrate the performance of the optimized myopic bounds over a range of transition probabilities.
 Consider a POMDP  with $\statedim=3$, $\actiondim = 2$, additive Gaussian noise model of Example 2, and   transition  matrices 
\begin{align}\nonumber
\tp_2 =\begin{pmatrix}
1 & 0 & 0\\
1-2\theta_1 & \theta_1 & \theta_1\\
1-2\theta_2 & \theta_2 & \theta_2
\end{pmatrix},\; \tp_1 = \tp_2^2.
\end{align}
It is straightforward to show that for all probabilities $\theta_1, \theta_2$ such that $\theta_1 + \theta_2 \le 1, \theta_2 \ge \theta_1$, 
the assumptions of Theorem \ref{th:theorem1} hold. (More generally choosing $\tp_1 = \tp_2^n$ for any positive integer $n$ satisfies the assumptions).
 We chose the cost vectors as  $\cost_1 = \left[1, 1.1, 1.2\right]^\prime$ and $\cost_2 = \left[1.2, 1.1, 1.1\right]^\prime$. 
Table \ref{tb:table1f} displays the worst case and best case values for performance metrics $(\text{vol}\left(\setoverlap\right),\percentloss_1,\percentloss_2)$ versus discount factor $\discount$
by sweeping over the entire range of $(\theta_1,\theta_2)$. The worst case performance is  denoted by  $\underline{\text{vol}}\left(\setoverlap\right)$,  $\underline{\percentloss}_1$,
 $\underline{\percentloss}_2$ and the best case  
by $\overline{\text{vol}}\left(\setoverlap\right)$, $\overline{\percentloss}_1$,
 $\overline{\percentloss}_2$.

%
%Parameters $\text{vol}\left(\setoverlap\right)$, $\percentloss_1$ and $\percentloss_2$, given in Table \ref{tb:table1f}, are obtained by averaging over 1000 randomly generated $\tp_1$ and $\tp_2$ for different values of discount factor $\discount$. $\percentloss_1$ is obtained by assuming $\belief_0 = \basisvec_3$ (myopic bounds overlap at $\basisvec_3$) and $\percentloss_2$ is obtained by sampling $\belief_0 \notin \setoverlap$.

\subsubsection*{Discussion}  In numerical examples, we found that
the percentage loss in optimality $\percentloss$ defined in (\ref{eq:percentloss})
depends on the following:
 whether or not the initial belief $\belief_0$ is in  $\setoverlap$,  $\text{vol}\left(\setoverlap\right)$, the trajectories of belief transitions and discount factor $\discount$.
 % and how long the sample path of posterior beliefs remain in $\setoverlap$.
 For $\belief_0 \notin \setoverlap$,
 one  expects $\percentloss$  to increase as $\discount$ decreases, and this is an observable trend in the above examples. 
 This is because as $\discount$ decreases, most of the value and loss in optimality is incurred in the near term.
  Symmetrically, if $\belief_0 \in \setoverlap$ one  expects $\percentloss$ to decrease as 
  $\discount$ decreases, because the  first few decisions determine most of the  incurred cost and 
  this relevant time period is  more likely to feature beliefs within $\setoverlap$.
   This is often the case but is not consistently so. The non-monotonicity of $\percentloss$ in these examples
derives from belief trajectories that can migrate out of and back into $\setoverlap$ as information accrues.
% In general,  $\percentloss$ depends on a combination of  whether $\belief_0 \in \setoverlap$,
%$\text{vol}\left(\setoverlap\right)$, the trajectories of belief transitions and discount factor $\discount$.

The electronic companion   to this paper  shows how the  results  be extended to problems with quadratic costs in the belief state. Such
problems arise in controlled sensing applications involving radars and sensor scheduling. Also further discussion  on  examples
that satisfy the assumptions of this paper is given.
 
%%%%%%%%%%

\section*{Appendix}
%\begin{definition}
Let $\belief_1, \belief_2 \in \Belief$, be any two belief states. Then, $\belief_1$ dominates $\belief_2$ with respect to MLR ordering, i.e. $\belief_1 \gr \belief_2$, if
$\belief_1(i)\belief_2(j) \le \belief_1(j)\belief_2(i), i < j,i,j \in \statespace$.
Also, $\belief_1$ dominates $\belief_2$ with respect to first order stochastic dominance, i.e. $\belief_1 \gs \belief_2$, if
$\sum_{i=q}^{\statedim}\belief_1(i) \ge \sum_{i=q}^{\statedim}\belief_2(i), q \in \statespace$.

The following lemma replaces Lemma 1.2(3) and Lemma 2.3(c) in \cite{Lov87}. The proof is straightforward and omitted.
\begin{lemma}\label{lm:norm_a1_ls_a2_gen}
Assume \textbf{(A\ref{it:norm_ac_inc})} holds. Then, for all $\belief \in \Belief, \filternorm(\belief,\action+1) \gs \filternorm(\belief,\action)$ where $\filternorm \left(\belief,\action \right) \equiv \left[\filternorm \left(\belief, 1, \action \right), \cdots,\filternorm \left(\belief,\obsdim,\action \right)\right]$.
\\ Assume \textbf{(A\ref{it:Copo_gen})} holds. Then, for all $\belief \in \Belief, \filter(\belief,\obs,\action+1) \gr \filter(\belief,\obs,\action)$.
\end{lemma}

{\textit{Proof of Theorem \ref{th:theorem1}}:} We show that
under \textbf{(A\ref{it:increasing_cost})}, \textbf{(A\ref{it:TP2_gen})}, \textbf{(A\ref{it:Copo_gen})} and \textbf{(A\ref{it:norm_ac_inc})} ,
$\optpolicy(\belief) \le \policyu(\belief)~\forall \belief \in \Belief$.
Let $\bvalue$ and $\bvalueaction$ denote the variables in Bellman's equation (\ref{eq:bellman}) when using costs $ \ucost_{\action}$ defined in (\ref{eq:new_cost_u}).
Then from \cite[Lemma 1.2.1]{Lov87} and \cite[Proposition 1]{Lov87}, $\bvalue(\filter(\belief,\obs,\action))$ is increasing in $\obs$. From Lemma \ref{lm:norm_a1_ls_a2_gen}, $\filternorm(\belief,\action+1) \gs \filternorm(\belief,\action)$. Therefore, % in \eqref{eq:proof_1}, (a) holds from the definition of $\gs$.
\begin{align}\label{eq:proof_1}
\begin{aligned}
\sum_{\obs \in \obspace}\bvalue(\filter(\belief,\obs,\action))\filternorm(\belief,\obs,\action)
\underset{(a)}{\le} &\sum_{\obs \in \obspace}\bvalue(\filter(\belief,\obs,\action))\filternorm(\belief,\obs,\action+1)
\underset{(b)}{\le} &\sum_{\obs \in \obspace}\bvalue(\filter(\belief,\obs,\action+1))\filternorm(\belief,\obs,\action+1)
\end{aligned}
\end{align}
 Inequality (b) holds since from Lemma \ref{lm:norm_a1_ls_a2_gen} and \cite[Proposition 1]{Lov87}, $\bvalue(\filter(\belief,\obs,\action+1)) \ge \bvalue(\filter(\belief,\obs,\action)) \forall \obs \in \obspace$.  The implication of \eqref{eq:proof_1} is that $\sum_{\obs \in \obspace}\bvalue(\filter(\belief,\obs,\action))\filternorm(\belief,\obs,\action)$ is increasing w.r.t $\action$ or equivalently,
\begin{align}\label{eq:proof_2}
\begin{aligned}
&\bvalueaction(\belief,\action) - \ucost_{\action}'\belief \le  \bvalueaction(\belief,\action+1) - \ucost_{\action+1}'\belief\\
&\Rightarrow \optpolicy(\belief) =\underset{\action \in \actionspace}{\argmin} ~\bvalueaction(\belief,\action) \le \underset{\action \in \actionspace}{\argmin} ~\ucost_{\action}'\belief = \policyu(\belief)
\end{aligned}
\end{align}
where the implication in \eqref{eq:proof_2} follows from \cite[Lemma 2.2]{Lov87}.
The proof that $\optpolicy(\belief) \ge \policyl(\belief)$ is similar and omitted.

{\textit{Proof of Theorem \ref{th:theorem_opt_bounds}}:}  The sufficient conditions in \eqref{eq:optimal_overlap} ensure that $\Belief^{\of^*}_{1} \supseteq \Belief^{\of}_{1} ~\forall \of \in \poly_{\of}$ and $\Belief^{\uf^*}_{2} \supseteq \Belief^{\uf}_{2} ~\forall \uf \in \poly_{\uf}$.
Indeed, to establish that $\text{vol}\left(\Belief^{\of^*}_{1}\right)$ $\ge$ $\text{vol}\Big(\Belief^{\of}_{1}\Big)$ $\forall \of \in \poly_{\of}$:
\begin{align}\label{eq:proof_th_opt_bounds_1}
\begin{aligned}
&\left(\tp_1 - \tp_2\right)\of^* \gc \left(\tp_1 - \tp_2\right)\of \quad \forall \of \in \poly_{\of}\\
\Rightarrow~& c_1 - c_2 - \discount\left(\tp_1 - \tp_2\right)\of^* \lc  c_1 - c_2 - \discount\left(\tp_1 - \tp_2\right)\of ~~\forall \of \in \poly_{\of}\\
\Rightarrow~& \Belief^{\of^*}_{1} \supseteq \Belief^{\of}_{1} ~\forall \of \in \poly_{\of}
\Rightarrow \text{vol}\left(\Belief^{\of^*}_{1}\right) \ge \text{vol}\Big(\Belief^{\of}_{1}\Big) \forall \of \in \poly_{\of}
\end{aligned}
\end{align}
So $\text{vol}\left(\Belief^{\of^*}_{1}\right) \ge \text{vol}\Big(\Belief^{\of}_{1}\Big) \forall \of \in \poly_{\of}$ and $\text{vol}\left(\Belief^{\uf^*}_{2}\right) \ge \text{vol}\Big(\Belief^{\uf}_{2}\Big) \forall \of \in \poly_{\of}$. Since $\setoverlap = \Belief^{\of^*}_{1} \cup \Belief^{\uf^*}_{2}$, the proof is complete.
\hrule
{\footnotesize
{\textit{Parameters of Example 1}}: For  the first example the parameters are defined as,
\begin{align} \nonumber
\begin{aligned}
c &=
 \begin{pmatrix}
    1.0000 &   1.5045   & 1.8341\\
    1.5002 &   1.0000   & 1.0000
 \end{pmatrix}^\prime,~
\tp_{2} =
 \begin{pmatrix}
1.0000  &  0.0000 &   0.0000\\
0.4677  &  0.4149 &   0.1174\\
0.3302  &  0.5220 &   0.1478
 \end{pmatrix}
,~\tp_1 = \tp_2^2\\
\oprob_{1} &=
 \begin{pmatrix}
 0.6373 &   0.3405 &   0.0222\\
 0.3118 &   0.6399 &   0.0483\\
 0.0422 &   0.8844 &   0.0734
 \end{pmatrix}
,~
 \oprob_2 =
 \begin{pmatrix}
 0.5927  &  0.3829 &   0.0244\\
0.4986  &  0.4625 &   0.0389\\
0.1395  &  0.79   &   0.0705
 \end{pmatrix}.
\end{aligned}
\end{align}

%{\textit{Parameters of Example 2}}: $\oprob_\action = \trid_{0.6} ~\forall \action \in \actionspace, \text{where}~\trid_\varepsilon~\text{is a tridiagonal matrix defined as}$
%    \begin{align} \nonumber
%    \begin{aligned}
%        ~\trid_\varepsilon = \left[\varepsilon_{ij}\right]_{\statedim \times \statedim}, \varepsilon_{ij} =\begin{cases} \varepsilon &i=j\\
%        1-\varepsilon &(i,j) = (1,2),(\statedim-1,\statedim)\\
%        \cfrac{1-\varepsilon}{2} &(i,j)=(i,i+1),(i,i-1),i\ne 1,\statedim\\
%        0 &\text{otherwise}\end{cases}
%          \end{aligned}
%    \end{align}
%\begin{align} \nonumber
%\begin{aligned}
%\tp_{1} =
% \begin{pmatrix}
%  1 & 0 & 0 \\
%  0.5901  &  0.3531  &  0.0568 \\
%  0.0179  &  0.8301   & 0.1521
% \end{pmatrix}
%, \tp_2 = \tp_1^2, \cost =
% \begin{pmatrix}
%  1.0457  &  1.7284   & 2.0487 \\
%  1.5457  &  1.0000  &  1.0000
% \end{pmatrix}^\p
%\end{aligned}
%\end{align}

{\textit{Parameters of Example 2}: For discrete observations $\oprob_\action = \oprob ~\forall \action \in \actionspace$,}
{\begin{align} \nonumber
    \begin{aligned}
    &\oprob =
     \begin{pmatrix}
0.0297 &   0.1334 &   0.1731 &   0.0482 &   0.1329 &   0.1095 &   0.0926 &   0.0348 &   0.1067 &   0.1391\\
0.0030 &   0.0271 &   0.0558 &   0.0228 &   0.0845 &   0.0923 &   0.1029 &   0.0511 &   0.2001 &   0.3604\\
0.0003 &   0.0054 &   0.0169 &   0.0094 &   0.0444 &   0.0599 &   0.0812 &   0.0487 &   0.2263 &   0.5075\\
     0 &   0.0011 &   0.0051 &   0.0038 &   0.0225 &   0.0368 &   0.0593 &   0.0418 &   0.2250 &   0.6046\\
     0 &   0.0002 &   0.0015 &   0.0015 &   0.0113 &   0.0223 &   0.0423 &   0.0345 &   0.2133 &   0.6731\\
     0 &        0 &   0.0005 &   0.0006 &   0.0056 &   0.0134 &   0.0298 &   0.0281 &   0.1977 &   0.7243\\
     0 &        0 &   0.0001 &   0.0002 &   0.0028 &   0.0081 &   0.0210 &   0.0227 &   0.1813 &   0.7638\\
     0 &        0 &        0 &   0.0001 &   0.0014 &   0.0048 &   0.0147 &   0.0183 &   0.1651 &   0.7956\\
     0 &        0 &        0 &        0 &   0.0007 &   0.0029 &   0.0103 &   0.0147 &   0.1497 &   0.8217\\
     0 &        0 &        0 &        0 &   0.0004 &   0.0017 &   0.0072 &   0.0118 &   0.1355 &   0.8434
     \end{pmatrix}
    \end{aligned}
\end{align}
\begin{align} \nonumber
    \begin{aligned}
    &\tp_{1} =
     \begin{pmatrix}
0.9496  &  0.0056 &   0.0056 &   0.0056&    0.0056 &   0.0056 &   0.0056 &   0.0056  &  0.0056 &   0.0056\\
0.9023  &  0.0081 &   0.0112 &   0.0112&    0.0112 &   0.0112 &   0.0112 &   0.0112  &  0.0112 &   0.0112\\
0.8574  &  0.0097 &   0.0166 &   0.0166&    0.0166 &   0.0166 &   0.0166 &   0.0166  &  0.0166 &   0.0167\\
0.8145  &  0.0109 &   0.0218 &   0.0218&    0.0218 &   0.0218 &   0.0218 &   0.0218  &  0.0218 &   0.0220\\
0.7737  &  0.0119 &   0.0268 &   0.0268&    0.0268 &   0.0268 &   0.0268 &   0.0268  &  0.0268 &   0.0268\\
0.7351  &  0.0126 &   0.0315 &   0.0315&    0.0315 &   0.0315 &   0.0315 &   0.0315  &  0.0315 &   0.0318\\
0.6981  &  0.0131 &   0.0361 &   0.0361&    0.0361 &   0.0361 &   0.0361 &   0.0361  &  0.0361 &   0.0361\\
0.6632  &  0.0136 &   0.0404 &   0.0404&    0.0404 &   0.0404 &   0.0404 &   0.0404  &  0.0404 &   0.0404\\
0.6301  &  0.0139 &   0.0445 &   0.0445&    0.0445 &   0.0445 &   0.0445 &   0.0445  &  0.0445 &   0.0445\\
0.5987  &  0.0141 &   0.0484 &   0.0484&    0.0484 &   0.0484 &   0.0484 &   0.0484  &  0.0484 &   0.0484
     \end{pmatrix}
    \end{aligned}
\end{align}
\begin{align} \nonumber
    \begin{aligned}
      \tp_{2} &=
          \begin{pmatrix}
0.5688 &   0.0143  &  0.0521  &  0.0521  &  0.0521 &   0.0521  &  0.0521 &   0.0521  &  0.0521 &   0.0522\\
0.5400 &   0.0144  &  0.0557  &  0.0557  &  0.0557 &   0.0557  &  0.0557 &   0.0557  &  0.0557 &   0.0557\\
0.5133 &   0.0145  &  0.0590  &  0.0590  &  0.0590 &   0.0590  &  0.0590 &   0.0590  &  0.0590 &   0.0592\\
0.4877 &   0.0145  &  0.0622  &  0.0622  &  0.0622 &   0.0622  &  0.0622 &   0.0622  &  0.0622 &   0.0624\\
0.4631 &   0.0145  &  0.0653  &  0.0653  &  0.0653 &   0.0653  &  0.0653 &   0.0653  &  0.0653 &   0.0653\\
0.4400 &   0.0144  &  0.0682  &  0.0682  &  0.0682 &   0.0682  &  0.0682 &   0.0682  &  0.0682 &   0.0682\\
0.4181 &   0.0144  &  0.0709  &  0.0709  &  0.0709 &   0.0709  &  0.0709 &   0.0709  &  0.0709 &   0.0712\\
0.3969 &   0.0143  &  0.0736  &  0.0736  &  0.0736 &   0.0736  &  0.0736 &   0.0736  &  0.0736 &   0.0736\\
0.3771 &   0.0141  &  0.0761  &  0.0761  &  0.0761 &   0.0761  &  0.0761 &   0.0761  &  0.0761 &   0.0761\\
0.3585 &   0.0140  &  0.0784  &  0.0784  &  0.0784 &   0.0784  &  0.0784 &   0.0784  &  0.0784 &   0.0787
        \end{pmatrix}
\\
    \end{aligned}
    \end{align}}
  {\begin{align} \nonumber
    \begin{aligned}
        \cost &=
          \begin{pmatrix}
    0.5986 &   0.5810  &  0.6116 &   0.6762  &  0.5664  &  0.6188  &  0.7107 &   0.4520 &   0.5986 &   0.7714\\
    0.6986 &   0.6727  &  0.7017 &   0.7649  &  0.6536  &  0.6005  &  0.6924 &   0.4324 &   0.5790 &   0.6714
          \end{pmatrix}^\prime
          \end{aligned}
    \end{align}}

{\textit{Parameters of Example 3}: $\oprob_\action = \trid_{0.7} ~\forall \action \in \actionspace$},
$\text{where}~\trid_\varepsilon~\text{is a tridiagonal matrix defined as}$
    \begin{align} \nonumber
    \begin{aligned}
        ~\trid_\varepsilon = \left[\varepsilon_{ij}\right]_{\statedim \times \statedim}, \varepsilon_{ij} =\begin{cases} \varepsilon &i=j\\
        1-\varepsilon &(i,j) = (1,2),(\statedim-1,\statedim)\\
        \cfrac{1-\varepsilon}{2} &(i,j)=(i,i+1),(i,i-1),i\ne 1,\statedim\\
        0 &\text{otherwise}\end{cases}
          \end{aligned}
    \end{align}
\begin{align} \nonumber
    \begin{aligned}
    &\tp_{1} =
     \begin{pmatrix}
   0.1851 &   0.1692 &   0.1630 &   0.1546 &   0.1324  &  0.0889 &   0.0546 &   0.0522\\
     0.1538 &   0.1531 &   0.1601 &   0.1580 &   0.1395  &  0.0994 &   0.0667 &   0.0694\\
     0.1307 &   0.1378 &   0.1489 &   0.1595 &   0.1472  &  0.1143 &   0.0769 &   0.0847\\
     0.1157 &   0.1307 &   0.1437 &   0.1591 &   0.1496  &  0.1199 &   0.0840 &   0.0973\\
     0.1053 &   0.1196 &   0.1388 &   0.1579 &   0.1520  &  0.1248 &   0.0888 &   0.1128\\
     0.0850 &   0.1056 &   0.1326 &   0.1618 &   0.1585  &  0.1348 &   0.0977 &   0.1240\\
     0.0707 &   0.0906 &   0.1217 &   0.1578 &   0.1629  &  0.1447 &   0.1078 &   0.1438\\
     0.0549 &   0.0757 &   0.1095 &   0.1502 &   0.1666  &  0.1576 &   0.1189 &   0.1666
     \end{pmatrix}
    \end{aligned}
\end{align}
\begin{align} \nonumber
    \begin{aligned}
      \tp_{2} &=
          \begin{pmatrix}
    0.0488 &   0.0696 &   0.1016 &   0.1413 &   0.1599 &   0.1614 &   0.1270  &  0.1904\\
     0.0413 &   0.0604 &   0.0882 &   0.1292 &   0.1503 &   0.1661 &   0.1425  &  0.2220\\
     0.0329 &   0.0482 &   0.0752 &   0.1195 &   0.1525 &   0.1694 &   0.1519  &  0.2504\\
     0.0248 &   0.0388 &   0.0649 &   0.1097 &   0.1503 &   0.1732 &   0.1643  &  0.2740\\
     0.0196 &   0.0309 &   0.0566 &   0.0985 &   0.1429 &   0.1805 &   0.1745  &  0.2965\\
     0.0158 &   0.0258 &   0.0517 &   0.0934 &   0.1392 &   0.1785 &   0.1794  &  0.3162\\
     0.0134 &   0.0221 &   0.0463 &   0.0844 &   0.1335 &   0.1714 &   0.1822  &  0.3467\\
     0.0110 &   0.0186 &   0.0406 &   0.0783 &   0.1246 &   0.1679 &   0.1899  &  0.3691
          \end{pmatrix}
    \end{aligned}
    \end{align}
    \begin{align} \nonumber
    \begin{aligned}
    &\tp_{3} =
     \begin{pmatrix}
   0.0077 &   0.0140 &   0.0337 &   0.0704  &  0.1178  &  0.1632  &  0.1983  &  0.3949\\
     0.0058 &   0.0117 &   0.0297 &   0.0659  &  0.1122  &  0.1568  &  0.1954  &  0.4225\\
     0.0041 &   0.0090 &   0.0244 &   0.0581  &  0.1011  &  0.1494  &  0.2013  &  0.4526\\
     0.0032 &   0.0076 &   0.0210 &   0.0515  &  0.0941  &  0.1400  &  0.2023  &  0.4803\\
     0.0022 &   0.0055 &   0.0165 &   0.0439  &  0.0865  &  0.1328  &  0.2006  &  0.5120\\
     0.0017 &   0.0044 &   0.0132 &   0.0362  &  0.0751  &  0.1264  &  0.2046  &  0.5384\\
     0.0012 &   0.0033 &   0.0106 &   0.0317  &  0.0702  &  0.1211  &  0.1977  &  0.5642\\
     0.0009 &   0.0025 &   0.0091 &   0.0273  &  0.0638  &  0.1134  &  0.2004  &  0.5826
     \end{pmatrix}
    \end{aligned}
\end{align}
\begin{align} \nonumber
    \begin{aligned}
      \tp_{4} &=
          \begin{pmatrix}
    0.0007 &   0.0020 &   0.0075 &   0.0244  &  0.0609  &  0.1104  &  0.2013  &  0.5928\\
     0.0005 &   0.0016 &   0.0063 &   0.0208  &  0.0527  &  0.1001  &  0.1991  &  0.6189\\
     0.0004 &   0.0013 &   0.0049 &   0.0177  &  0.0468  &  0.0923  &  0.1981  &  0.6385\\
     0.0003 &   0.0009 &   0.0038 &   0.0149  &  0.0407  &  0.0854  &  0.2010  &  0.6530\\
     0.0002 &   0.0007 &   0.0031 &   0.0123  &  0.0346  &  0.0781  &  0.2022  &  0.6688\\
     0.0001 &   0.0005 &   0.0023 &   0.0100  &  0.0303  &  0.0713  &  0.1980  &  0.6875\\
     0.0001 &   0.0004 &   0.0019 &   0.0083  &  0.0266  &  0.0683  &  0.1935  &  0.7009\\
     0.0001 &   0.0003 &   0.0014 &   0.0069  &  0.0240  &  0.0651  &  0.1878  &  0.7144
          \end{pmatrix}
    \end{aligned}
    \end{align}
    \begin{align} \nonumber
    \begin{aligned}
    &\tp_{5} =
     \begin{pmatrix}
        0.0000 &   0.0002 &   0.0010 &   0.0054  &  0.0204  &  0.0590  &  0.1772  &  0.7368\\
     0.0000 &   0.0001 &   0.0008 &   0.0041  &  0.0168  &  0.0515  &  0.1663  &  0.7604\\
     0.0000 &   0.0001 &   0.0006 &   0.0038  &  0.0156  &  0.0480  &  0.1596  &  0.7723\\
     0.0000 &   0.0001 &   0.0005 &   0.0032  &  0.0139  &  0.0450  &  0.1603  &  0.777\\
     0.0000 &   0.0001 &   0.0004 &   0.0028  &  0.0124  &  0.0418  &  0.1590  &  0.7835\\
     0.0000 &   0.0001 &   0.0003 &   0.0023  &  0.0106  &  0.0389  &  0.1547  &  0.7931\\
     0.0000 &   0.0000 &   0.0003 &   0.0018  &  0.0090  &  0.0351  &  0.1450  &  0.8088\\
     0.0000 &   0.0000 &   0.0002 &   0.0015  &  0.0080  &  0.0325  &  0.1386  &  0.8192
     \end{pmatrix}
    \end{aligned}
\end{align}
\begin{align} \nonumber
    \begin{aligned}
      \tp_{6} &=
          \begin{pmatrix}
    0.0000 &   0.0000 &   0.0001 &   0.0012  &  0.0067  &  0.0296  &  0.1331  &  0.8293\\
     0.0000 &   0.0000 &   0.0001 &   0.0010  &  0.0059  &  0.0275  &  0.1238  &  0.8417\\
     0.0000 &   0.0000 &   0.0001 &   0.0009  &  0.0056  &  0.0272  &  0.1238  &  0.8424\\
     0.0000 &   0.0000 &   0.0001 &   0.0009  &  0.0053  &  0.0269  &  0.1234  &  0.8434\\
     0.0000 &   0.0000 &   0.0001 &   0.0006  &  0.0043  &  0.0237  &  0.1189  &  0.8524\\
     0.0000 &   0.0000 &   0.0001 &   0.0005  &  0.0038  &  0.0215  &  0.1129  &  0.8612\\
     0.0000 &   0.0000 &   0.0000 &   0.0004  &  0.0032  &  0.0191  &  0.1094  &  0.8679\\
     0.0000 &   0.0000 &   0.0000 &   0.0003  &  0.0025  &  0.0161  &  0.1011  &  0.8800
          \end{pmatrix}
    \end{aligned}
    \end{align}
    \begin{align} \nonumber
    \begin{aligned}
    &\tp_{7} =
     \begin{pmatrix}
   0.0000 &   0.0000 &   0.0000 &   0.0003  &  0.0022  &  0.0143  &  0.0938  &  0.8894\\
     0.0000 &   0.0000 &   0.0000 &   0.0002  &  0.0019  &  0.0136  &  0.0901  &  0.8942\\
     0.0000 &   0.0000 &   0.0000 &   0.0002  &  0.0017  &  0.0126  &  0.0849  &  0.9006\\
     0.0000 &   0.0000 &   0.0000 &   0.0002  &  0.0015  &  0.0118  &  0.0819  &  0.9046\\
     0.0000 &   0.0000 &   0.0000 &   0.0001  &  0.0013  &  0.0108  &  0.0754  &  0.9124\\
     0.0000 &   0.0000 &   0.0000 &   0.0001  &  0.0011  &  0.0098  &  0.0714  &  0.9176\\
     0.0000 &   0.0000 &   0.0000 &   0.0001  &  0.0010  &  0.0090  &  0.0713  &  0.9186\\
     0.0000 &   0.0000 &   0.0000 &   0.0001  &  0.0009  &  0.0084  &  0.0675  &  0.9231
     \end{pmatrix}
    \end{aligned}
\end{align}
\begin{align} \nonumber
    \begin{aligned}
      \tp_{8} &=
          \begin{pmatrix}
    0.0000 &   0.0000 &   0.0000 &   0.0001  &  0.0008  &  0.0078  &  0.0665  &  0.9248\\
     0.0000 &   0.0000 &   0.0000 &   0.0000  &  0.0007  &  0.0068  &  0.0626  &  0.9299\\
     0.0000 &   0.0000 &   0.0000 &   0.0000  &  0.0006  &  0.0061  &  0.0581  &  0.9352\\
     0.0000 &   0.0000 &   0.0000 &   0.0000  &  0.0005  &  0.0057  &  0.0561  &  0.9377\\
     0.0000 &   0.0000 &   0.0000 &   0.0000  &  0.0005  &  0.0053  &  0.0558  &  0.9384\\
     0.0000 &   0.0000 &   0.0000 &   0.0000  &  0.0004  &  0.0051  &  0.0558  &  0.9387\\
     0.0000 &   0.0000 &   0.0000 &   0.0000  &  0.0004  &  0.0045  &  0.0522  &  0.9429\\
     0.0000 &   0.0000 &   0.0000 &   0.0000  &  0.0003  &  0.0040  &  0.0505  &  0.9452
          \end{pmatrix}
    \end{aligned}
    \end{align}
 \begin{align} \nonumber
    \begin{aligned}
      \cost &=
          \begin{pmatrix}
     1.0000 &   2.2486 &   4.1862 &   6.9509 &  11.2709 &  15.9589 &  21.4617 &  27.6965 \\
    31.3230 &   8.8185 &   9.6669 &  11.4094 &  14.2352 &  17.8532 &  22.3155 &  27.5353 \\
    50.0039 &  26.3162 &  14.6326 &  15.3534 &  17.1427 &  19.7455 &  23.1064 &  27.3025 \\
    65.0359 &  40.2025 &  27.5380 &  19.5840 &  20.3017 &  21.8682 &  24.2022 &  27.4108 \\
    79.1544 &  53.1922 &  39.5408 &  30.5670 &  23.3697 &  23.9185 &  25.1941 &  27.4021 \\
    90.7494 &  63.6983 &  48.6593 &  38.6848 &  30.4868 &  25.7601 &  26.0012 &  27.1867 \\
    99.1985 &  71.1173 &  55.0183 &  44.0069 &  34.7860 &  29.0205 &  26.9721 &  27.1546 \\
  106.3851  & 77.2019  & 60.0885  & 47.8917  & 37.6330  & 30.8279  & 27.7274  & 26.4338
          \end{pmatrix}
    \end{aligned}
    \end{align}}

    \ACKNOWLEDGMENT{This research was supported by  an insight grant from the  Social Sciences and Humanities Research Council of Canada,
    a discovery grant from the National Science and Engineering Research Council of Canada, and the Canada Research Chairs Program.}

%\bibliographystyle{ormsv080}
%\bibliography{$HOME/styles/bib/vkm}

\end{document}